\documentclass[12pt]{article}
\usepackage[latin1]{inputenc}
\usepackage{amssymb,amsmath,amsfonts}
\usepackage{xcolor}
\newtheorem{theorem}{Theorem}
\newtheorem{proposition}[theorem]{Proposition}

\newtheorem{example}[theorem]{Example}

\newcommand{\R}{\mathbb{R}}
\newcommand{\Sp}{\mathbb{S}}

\renewcommand{\L}{\mathbb{L}}

\newcommand{\spa}{\mbox{span\,}}

\newcommand{\Hy}{\mathbb{H}}
\newcommand{\rank}{\mbox{rank }}

\newcommand{\po}{{\hspace*{-1ex}}{\bf .  }}

\def\lp{{\langle\!\langle}}\vspace{2ex}
\def\rp{{\rangle\!\rangle}}
\def\<{{\langle}}
\def\>{{\rangle}}
\def\Sal{{\cal S}}

\def\a{\alpha}
\def\be{\begin{equation} }
\def\ee{\end{equation} }

\def\proof{\noindent\emph{Proof: }}
\def\qed{\ifhmode\unskip\nobreak\fi\ifmmode\ifinner
\else\hskip5 pt \fi\fi\hbox{\hskip5 pt \vrule width4 pt
height6 pt  depth1.5 pt \hskip 1pt }}
\newcommand\blfootnote[1]{
\begingroup
\renewcommand\thefootnote{}\footnote{#1}
\addtocounter{footnote}{-1}
\endgroup
}
\begin{document}

\title{Conformal Kaehler submanifolds}
\author{L. J. Al\'ias\footnote{Corresponding author. e-mail: ljalias@um.es},  \, S. Chion and M. Dajczer}
\date{}
\maketitle

\begin{abstract}  This paper presents two results in the realm 
of conformal Kaehler submanifolds. These are conformal immersions 
of Kaehler manifolds into the standard flat Euclidean space. 
The proofs are obtained by making a rather strong use of several 
facts and techniques developed in \cite{CD} for the study of 
isometric immersions of Kaehler manifolds into the standard 
hyperbolic space. 
\end{abstract}
\blfootnote{\textup{2020} \textit{Mathematics Subject Classification}:
53B25, 53C42.}

Let $f\colon M^{2n}\to\R^{2n+p}$ denote a \emph{conformal Kaehler 
submanifold}. Thus $(M^{2n},J)$ is a Kaehler manifold of complex 
dimension $n\geq 2$ and $f$  a conformal immersion into Euclidean 
space that lies in codimension $p$. Thus, there is a positive 
function $\lambda \in C^\infty(M)$ such that 
the Kaehler metric and the one induced by $f$ 
relate by $\<\,,\,\>_f=\lambda^2\<\,,\,\>_{M^{2n}}$. 

Conformal Kaehler submanifolds laying in the low codimensions $p=1$ 
and $p=2$ have already been considered in \cite{CCD}. In 
this paper, we are interested in higher codimensions 
although not too large in comparison to the dimension of the 
manifold. 

Our first result, provides a necessary condition for the existence
of a conformal immersion in codimension at most $n-3$ in 
terms of the sectional curvature of the Kaehler manifold.

\begin{theorem}\po\label{one}
Let $f\colon M^{2n}\to\R^{2n+p}$, $p\leq n-3$, be a conformal
Kaehler submanifold. Then at any $x\in M^{2n}$ there is a complex 
vector subspace $V^{2m}\subset T_xM$ with $m\geq n-p$ such that 
the sectional curvature of $M^{2n}$ satisfies $K_M(S,JS)(x)\leq 0$ 
for any $S\in V^{2m}$.
\end{theorem}

Notice that the conclusion of Theorem \ref{one} remains valid if 
the Euclidean ambient space $\R^{2n+p}$ is replaced by any locally 
conformally flat manifold of the same dimension.
\vspace{1ex}

Our second result characterizes a submanifold, in terms of the 
degree of positiveness of the sectional curvature, as being  
locally the Example \ref{example} presented below.
\vspace{1ex}

The light-cone $\mathbb{V}^{m+1}\subset\L^{m+2}$ of the standard 
flat Lorentzian space is any one of the two connected components 
of the set of all light-like vectors, namely,  
$$
\{v\in\L^{m+2}:\<v,v\>=0,\;v\neq 0\}
$$ 
endowed with the (degenerate) induced metric. The Euclidean 
space $\R^m$ can be realized as an umbilic hypersurface of 
$\mathbb{V}^{m+1}$ as follows: 
Take vectors $v,w\in\mathbb{V}^{m+1}$ such that 
$\<v,w\>=1$ and a linear isometry $C\colon\R^m\to\{v,w\}^\perp$. 
Let $\psi\colon\R^m\to\mathbb{V}^{m+1}\subset\L^{m+2}$ be 
defined by
\be\label{psi}
\psi(x)=v+Cx-\frac{1}{2}\|x\|^2w.
\ee
Then $\psi$ is an isometric embedding of $\R^m$ as an umbilical 
hypersurface in the light cone which is  the intersection of 
$\mathbb{V}^{m+1}$ with an affine hyperplane. Namely, we have that
$$
\psi(\R^m)=\{y\in\mathbb{V}^{m+1}\colon\<y,w\>=1\}.
$$
The normal bundle of $\psi$ is $N_\psi\R^m=\spa\{\psi,w\}$ 
and its second fundamental form is 
$$
\a^\psi(X,Y)=-\<X,Y\>_{\R^m}w.
$$

Proposition $9.9$ in \cite{DT} gives an elementary correspondence 
between the conformal immersions in Euclidean space and the isometric 
immersions into the light cone which goes as follows: Associated to 
a given conformal immersion $f\colon M^m\to\R^{m+p}$ with 
conformal factor $\lambda\in C^\infty(M)$ there is the associated
isometric immersion defined by
$$
F=\frac{1}{\lambda}\psi\circ f\colon M^m
\to\mathbb{V}^{m+p+1}\subset\L^{m+p+2}.
$$
Conversely, any isometric immersion 
$F\colon M^m\to\mathbb{V}^{m+p+1}\setminus\R w\subset\L^{m+p+2}$ 
gives rise to an associated conformal immersion $f\colon M^m\to\R^{m+p}$ 
given by $\psi\circ f=\pi\circ F$ with conformal factor $1/\<F,w\>$. 
Here $\pi\colon\mathbb{V}^{m+p+1}\setminus\R w\to\R^{m+p}$ 
is the projection $\pi(x)=x/\<x,w\>$.

\begin{example}\po\label{example}
{\em  Let the Kaehler manifold $M^{2n}$ be the Riemannian 
product of one hyperbolic plane and a set of
two-dimensional round spheres such that
$$
M^{2n}=\Hy^2_c\times\Sp^2_{c_2}\times\cdots\times\Sp^2_{c_n}
\;\;\mbox{with}\;\;1/c_2+\cdots+1/c_n=-1/c.
$$
If $f_1$ is the inclusion $\Hy^2_c\subset\L^3$ and 
$f_2\colon\Sp^2_{c_2}\times\cdots\times\Sp^2_{c_n}\to\Sp^{3n-4}_c
\subset\R^{3n-3}$ is the product of umbilical spheres then the 
map $\psi^{-1}\circ (f_1\times f_2)\colon M^{2n}\to\R^{3n-2}$ is a 
conformal Kaehler submanifold.
}\end{example}

\begin{theorem}\po\label{three}
Let $f\colon M^{2n}\to\R^{2n+p}$, $2\leq p\leq n-2$, be a 
connected conformal Kaehler submanifold. Assume that at a 
point $x_0\in M^{2n}$ there is a complex tangent vector 
subspace $V^{2m}\subset T_{x_0}M$ with $m\geq p+1$ such that the 
sectional curvature of $M^{2n}$ satisfies $K_M(S,JS)(x)>0$ for 
any $0\neq S\in V^{2m}$. Then $p=n-2$ and $f(M)$ is an open 
subset of the submanifold given by Example \ref{example}. 
\end{theorem}

\section{The proofs}

Let $V^{2n}$ and $\L^p$, $p\geq 2$, be real vector spaces 
such that there is $J\in Aut(V)$ which satisfies $J^2=-I$ 
and $\L^p$ is endowed with a Lorentzian inner product $\<\,,\,\>$. 
Then let $W^{p,p}=\L^p\oplus\L^p$ be endowed 
with the inner product of signature $(p,p)$ defined by
$$
\lp(\xi,\bar\xi),(\eta,\bar\eta)\rp
=\<\xi,\eta\>-\<\bar\xi,\bar\eta\>.
$$
A vector subspace $L\subset W^{p,p}$ is 
called \emph{degenerate} if $L\cap L^\perp\neq 0$.
\vspace{1ex}

Let $\a\colon V^{2n}\times V^{2n}\to\L^p$ be 
a symmetric bilinear form and 
$\beta\colon V^{2n}\times V^{2n}\to W^{p,p}$ the associated 
bilinear form given by
$$
\beta(X,Y)=(\a(X,Y)+\a(JX,JY),\a(X,JY)-\a(JX,Y)).
$$ 
We have that if $\beta(X,Y)=(\xi,\eta)$ then
\be\label{symmetries}
\beta(X,JY)=(\eta,-\xi)\;\,\mbox{and}\;\,\beta(Y,X)=(\xi,-\eta).
\ee
We denote the vector subspace of $W^{p,p}$ generated by 
$\beta$ by
$$
\mathcal{S}(\beta)=\spa\{\beta(X,Y)\colon X,Y\in V^{2n}\}
$$
and say that $\beta$ is \emph{surjective} if 
$\mathcal{S}(\beta)=W^{p,p}$. The (right) kernel $\beta$ 
is defined by
$$
\mathcal{N}(\beta)=\{Y\in V^{2n}\colon\beta(X,Y)=0
\;\,\mbox{for all}\;\,X\in V^{2n}\}.
$$

A vector $X\in V^{2n}$ is called a (left) \emph{regular element}
of $\beta$ if $\dim B_X(V)=r$ where
$r=\max\{\dim B_X(V)\colon X\in V\}$ and $B_X\colon V\to W^{p,p}$ 
is the linear transformation defined by $B_XY=\beta(X,Y)$.
The set $RE(\beta)$ of regular elements of
$\beta$ is easily seen to be an open dense subset
of $V^{2n}$, for instance see  Proposition $4.4$ in \cite{DT}.

It is said that $\beta$ is \emph{flat} if it satisfies that
$$
\lp\beta(X,Y),\beta(Z,T)\rp-\lp\beta(X,T),\beta(Z,Y)\rp=0
\;\,\mbox{for all}\;\,X,Z,Y,T\in V^{2n}.
$$
If $\beta$ is flat and $X\in RE(\beta)$ we have
from Proposition $4.6$ in \cite{DT} that
\be\label{moore}
\Sal(\beta|_{V\times\ker\beta_X})\subset B_X(V)\cap(B_X(V))^\perp.
\ee

\begin{proposition}\label{mainlemma}\po Let 
$\beta\colon V^{2n}\times V^{2n}\to W^{p,p}$, $p\leq n$, 
be flat and surjective. Then 
$$
\dim\mathcal{N}(\beta)\geq 2n-2p.
$$
\end{proposition}

\proof  This is condition $(9)$ in Proposition $11$ 
of \cite{CD}.
\vspace{2ex}\qed

Let $V^{2n}$ be endowed with a positive definite inner product 
$(\,,\,)$ with respect to which $J\in Aut(V)$ is an isometry. 
Assume that there is a light-like vector $w\in\L^p$ such that
\be\label{shapeid}
\<\a(X,Y),w\>=-(X,Y)\;\,\mbox{for any}\;\,X,Y\in V^{2n}.
\ee

Let $U_0^s\subset\L^p$ be the $s$-dimensional vector subspace 
given by 
$$
U_0^s=\pi_1(\Sal(\beta))=\spa\{\alpha(X,Y)
+\alpha(X,JY)\colon X,Y\in V^{2n}\},
$$
where $\pi_1\colon W^{p,p}\to\L^p$ 
denotes the projection onto the first component of $W^{p,p}$.  
From Proposition $9$ in \cite{CD} we know that
\be\label{ebetaimage}
\Sal(\beta)=U_0^s\oplus U_0^s. 
\ee
In addition, if $\Sal(\beta)$ is a degenerate vector subspace 
then $1\leq s\leq p-1$ and there is a light-like vector 
$v\in U_0^s$ such that
\be\label{uimage}
\Sal(\beta)\cap\Sal(\beta)^\perp=\spa\{v\}\oplus\spa\{v\}.
\ee

\begin{proposition}\po\label{decompor}
Let the bilinear form $\beta\colon V^{2n}\times V^{2n}\to W^{p,p}$ 
be flat and the vector subspace $\Sal(\beta)$ degenerate.
Then $L=\spa\{v,w\}\subset\L^p$ is a Lorentzian plane. 
Moreover, choosing $v$ such that $\<v,w\>=-1$ and setting
$\beta_1=\pi_{L^\perp\times L^\perp}\circ\beta$, we have
\be\label{betadecomp}
\beta(X,Y)=\beta_1(X,Y)+2((X,Y)v,(X,JY)v)
\;\,\mbox{for any}\;\, X,Y\in V^{2n}.
\ee 
Furthermore, if $s\leq n$ then
\be\label{estimnucleo}
\dim\mathcal{N}(\beta_1)\geq 2n-2s+2.
\ee
\end{proposition}

\proof We obtain from \eqref{uimage} that
\be\label{vdeg}
0=\lp\beta(X,Y),(v,0)\rp=\<\a(X,Y)+\a(JX,JY),v\>
\;\,\mbox{for any}\;\,X,Y\in V^{2n}.
\ee 
From \eqref{shapeid} and the fact that $J$ is an isometry with 
respect to $(\,,\,)$, we also have that
\be\label{luis1} 
\lp\beta(X,Y),(w,0)\rp=\<\a(X,Y)+\a(JX,JY),w\>=-2(X,Y)
\;\,\mbox{for any}\;\,X,Y\in V^{2n}.
\ee 
In particular, we have
$\lp\beta(X,X),(w,0)\rp<0$ for any $0\neq X\in V^{2n}$, which 
jointly with \eqref{vdeg} implies that $v$ and $w$ 
are linearly independent and thus span a Lorentzian plane.

Since $w$ is light-like and $v$ satisfies $\<v,w\>=-1$, we have
\begin{align*}
\a(X,Y)+\a(JX,JY)&=\a_{L^\perp}(X,Y)+\a_{L^\perp}(JX,JY)
-\<\a(X,Y)+\a(JX,JY),v\>w\\
&\quad-\<\a(X,Y)+\a(JX,JY),w\>v,
\end{align*}
where $\a_{L^\perp}$ denotes the $L^\perp$-component of $\a$. 
Then \eqref{vdeg} and \eqref{luis1} yield
$$
\a(X,Y)+\a(JX,JY)=\a_{L^\perp}(X,Y)+\a_{L^\perp}(JX,JY)+2(X,Y)v
$$ 
and
$$
\a(X,JY)-\a(JX,Y)=\a_{L^\perp}(X,JY)-\a_{L^\perp}(JX,Y)+2(X,JY)v,
$$
from which we obtain \eqref{betadecomp}. 

We have from \eqref{ebetaimage} and \eqref{uimage} that 
$w\notin U_0^s+L^\perp$. Hence $\dim(U_0^s+L^\perp)= p-1$.  
It then follows from
$$
p-1=\dim(U_0^s+L^\perp)=\dim U_0^s+\dim L^\perp-\dim U_0^s\cap L^\perp
$$
that $U_1=U_0^s\cap L^\perp$ satisfies 
\be\label{dim}
\dim U_1=s-1
\ee
and we have from \eqref{ebetaimage}, \eqref{uimage} and \eqref{betadecomp}
that $\Sal(\beta_1)=U_1^{s-1}\oplus U_1^{s-1}$.

From \eqref{betadecomp} we obtain that 
$$
\lp\beta(X,Y),\beta(Z,T)\rp=\lp\beta_1(X,Y),\beta_1(Z,T)\rp
\;\,\mbox{for any}\;\,X,Y,Z,T\in V^{2n},
$$ 
and hence also the  bilinear form
$\beta_1\colon V^{2n}\times V^{2n}\to L^\perp\oplus L^\perp$
is flat. Let $X\in RE(\beta_1)$ and set 
$N_1(X)=\ker B_{1X}$ where $B_{1X}Y=\beta_1(X,Y)$. To
obtain \eqref{estimnucleo} it suffices to show that 
$N_1(X)=\mathcal{N}(\beta_1)$ since then 
$\dim\mathcal{N}(\beta_1)
=\dim N_1(X)\geq 2n-2\dim U_1= 2n-2s+2$.

If $\beta_1(Y,Z)=(\xi,\eta)$ then by \eqref{symmetries} and \eqref{betadecomp} 
we have $\beta_1(Z,Y)=(\xi,-\eta)$.
If $Y,Z\in N_1(X)$ it follows from \eqref{moore} that
$$
0=\lp\beta_1(Y,Z),\beta_1(Z,Y)\rp
=\lp(\xi,\eta),(\xi,-\eta)\rp=\|\xi\|^2+\|\eta\|^2.
$$
Hence $\beta_1|_{N_1(X)\times N_1(X)}=0$ since the inner product 
induced on $U_1^{s-1}$ is positive definite. Now let 
$\beta_1(Y,Z)=(\delta,\zeta)$ where $Y\in V^{2n}$ and 
$Z\in N_1(X)$. Then the flatness of $\beta_1$ yields
$$
0=\lp\beta_1(Y,Z),\beta_1(Z,Y)\rp
=\lp(\delta,\zeta),(\delta,-\zeta)\rp
=\|\delta\|^2+\|\zeta\|^2
$$
and therefore $\beta_1|_{V\times N_1(X)}=0$.\qed

\begin{proposition}\po\label{umbilical}
Let the bilinear form $\beta\colon V^{2n}\times V^{2n}\to W^{p,p}$
be flat and satisfy 
\be\label{product}
\lp\beta(X,Y),\gamma(Z,T)\rp=\lp\beta(X,T),\gamma(Z,Y)\rp
\;\,\mbox{for any}\;\,X,Y,Z,T\in V^{2n}
\ee 
where $\gamma\colon V^{2n}\times V^{2n}\to W^{p,p}$
is the bilinear form defined by
$$
\gamma(X,Y)=(\a(X,Y),\a(X,JY))\;\,\mbox{for any}\;\,X,Y\in V^{2n}.
$$
If the vector subspace $\Sal(\beta)$ is degenerate and 
$s\leq n-1$ then there is a $J$-invariant vector subspace 
$P^{2m}\subset V^{2n}$, $m\geq n-s+1$, such that 
$$
\<\a(S,S),\a(JS,JS)\>-\|\a(S,JS)\|^2\leq 0\;\,\mbox{for any}\;\,S\in P^{2m}.
$$ 
\end{proposition}

\proof Let $v\in U_0^s$ be given by \eqref{uimage}. 
We claim that
\be\label{umbilic}
\<\a(X,Y),v\>=0\;\;\mbox{for any}\;\; X,Y\in V^{2n}.
\ee
Since $s\leq n-1$ then \eqref{estimnucleo} gives
$\dim\mathcal{N}(\beta_1)\geq 4$. Hence \eqref{betadecomp} 
yields $\beta(S,S)=2((S,S)v,0)$ for any $S\in\mathcal{N}(\beta_1)$. 
Thus
\be\label{product2}
\lp\gamma(X,Y),\beta(S,S)\rp=2\<\a(X,Y),v\>
\ee
for any $S\in\mathcal{N}(\beta_1)$ of unit length. 
On the other hand, we obtain from \eqref{symmetries}  and
\eqref{betadecomp} that $\beta(S,Y)=\beta(Y,S)=0$ for any
$S\in\mathcal{N}(\beta_1)$ and $Y\in\{S,JS\}^\perp$. 
Then \eqref{product} and \eqref{product2} give $\<\a(X,Y),v\>=0$
for any $X\in V^{2n}$ and $Y\in\{S,JS\}^\perp$ where
$S\in\mathcal{N}(\beta_1)$. Now that $\dim\mathcal{N}(\beta_1)\geq 4$
yields the claim. 

Choosing $v\in U_0^s$ as in Proposition \ref{decompor} it
follows from \eqref{shapeid} and \eqref{umbilic} that
\be\label{alphapar}
\a(X,Y)=\a_{L^\perp}(X,Y)+(X,Y)v\;\,\mbox{for any}\;\,X,Y\in V^{2n}.
\ee 
Then we obtain from \eqref{alphapar} that
$$
\gamma(X,Y)=
(\a_{L^\perp}(X,Y)+(X,Y)v,\a_{L^\perp}(X,JY)+(X,JY)v)
\;\,\mbox{for any}\;\,X,Y\in V^{2n}.
$$
Set $P^{2m}=\mathcal{N}(\beta_1)$ where
$2m=\dim\mathcal{N}(\beta_1)\geq 2n-2s+2$ by \eqref{estimnucleo}.
From \eqref{betadecomp} we have
$\beta(Z,S)=2((Z,S)v,(Z,JS)v)$ for any $S\in P^{2m}$
and $Z\in V^{2n}$. Then \eqref{product} gives
$$
\lp\gamma(X,S),\beta(Z,Y)\rp=\lp\gamma(X,Y),\beta(Z,S)\rp=0
$$
for any $S\in P^{2m}$ and $X,Y,Z\in V^{2n}$.
Hence  $\Sal(\gamma|_{V\times P})$
and $\Sal(\beta)$ are orthogonal vector subspaces.  From \eqref{dim}
we have
\be\label{deco}
U_0^s=U_1^{s-1}\oplus\spa\{v\}\;\,\mbox{where}\;\,U_1^{s-1}=U_0^s\cap L^\perp.
\ee
Then by \eqref{ebetaimage} the vector subspaces
$\Sal(\gamma|_{V\times P})$ and $U_1^{s-1}\oplus U_1^{s-1}$ are
orthogonal, and thus
$$
\<\a(X,S),\xi\>=\lp\gamma(X,S),(\xi,0)\rp=0
$$
for any $X\in V^{2n}$, $S\in P^{2m}$ and $\xi\in U^{s-1}_1$.
Since $U_1^{s-1}\subset L^\perp$ then
\be\label{nucleoa}
\a_{U_1}(X,S)=0\;\,\mbox{for any}\;\,X\in V^{2n}\;\mbox{and}\;S\in P^{2m}.
\ee

Let $\L^p=U_1^{s-1}\oplus U_2^{p-s-1}\oplus L$
be an orthogonal decomposition.
Then \eqref{ebetaimage} and \eqref{deco} give
$$
\<\a(X,Y)+\a(JX,JY),\xi_2\>
=\lp\beta(X,Y),(\xi_2,0)\rp=0
$$
for any $X,Y\in V^{2n}$ and $\xi_2\in U_2^{p-s-1}$.  Thus
\be\label{pluri}
\a_{U_2}(X,Y)=-\a_{U_2}(JX,JY)\;\,\mbox{for any}\;\,X,Y\in V^{2n}.
\ee

Having $U_2^{p-s-1}$ a positive definite induced inner product,
we obtain from \eqref{alphapar}, \eqref{nucleoa} and \eqref{pluri} 
that
$$
\<\a(S,S),\a(JS,JS)\>-\|\a(S,JS)\|^2
=-\|\a_{U_2}(S,S)\|^2-\|\a_{U_2}(S,JS)\|^2\leq 0
$$
for any $S\in P^{2m}$.\vspace{2ex}\qed 

Given a conformal immersion $f\colon M^{2n}\to\R^{2n+p}$ with 
conformal factor $\lambda\in C^\infty(M)$ we have the associated
isometric immersion $F=\frac{1}{\lambda}\psi\circ f\colon M^{2n}
\to\mathbb{V}^{2n+p+1}\subset\L^{2n+p+2}$ where $\psi$ is given
by \eqref{psi}. Differentiating $\<F,F\>=0$ once gives 
$F\in\Gamma(N_FM)$ and twice yields that the second fundamental 
form $\alpha^F\colon TM\times TM\to N_FM$ of $F$ satisfies
\be\label{condition}
\<\alpha^F(X,Y),F\>=-\<X,Y\>\;\,\mbox{for any}\;\,X,Y\in\mathfrak{X}(M).
\ee
Since $\psi_*N_fM\subset N_FM$ the normal bundle of $F$ decomposes as 
$N_FM=\psi_*N_fM\oplus L^2\equiv\mathbb{L}^{p+2}$ where $L^2$ is the Lorentzian plane 
subbundle orthogonal to $\psi_*N_fM$ such that $F\in\Gamma(L^2)$.

Let the bilinear forms  
$\gamma,\beta\colon T_xM\times T_xM\to N_FM(x)\oplus N_FM(x)$
be defined by 
$$
\gamma(X,Y)=(\a^F(X,Y),\a^F(X,JY))
$$ 
and
$$
\beta(X,Y)=(\a^F(X,Y)+\a^F(JX,JY),\a^F(X,JY)-\a^F(JX,Y)).
$$ 
\begin{proposition}\po\label{flatforms}
Let $N_FM(x)\oplus N_FM(x)$ be endowed with the inner product
defined by
$$
\lp(\xi,\bar\xi),(\eta,\bar\eta)\rp
=\<\xi,\eta\>-\<\bar\xi,\bar\eta\>.
$$
Then the bilinear form $\beta$ is flat and 
$$
\lp\beta(X,Y),\gamma(Z,T)\rp=\lp\beta(X,T),\gamma(Z,Y)\rp
\;\;\mbox{for any}\;\;X,Y,Z,T\in T_xM.
$$
\end{proposition}

\proof  The proof is straightforward using that 
$\beta(X,JY)=-\beta(JX,Y)$, that the curvature tensor satisfies 
$R(X,Y)JZ=JR(X,Y)Z$ for any $X,Y,Z\in T_xM$ and the Gauss equation 
for $f$; for details see the proof of Proposition $16$ in \cite{CD}
\vspace{2ex}\qed

\noindent{\em Proof of Theorem \ref{one}:} 
It suffices to show that the vector subspace $\Sal(\beta)$ is 
degenerate since then the proof follows from the Gauss equation 
jointly with Proposition \ref{umbilical} and Proposition \ref{flatforms}. 
If $\Sal(\beta)$ is not degenerate, and since we have the result 
given by Proposition \ref{flatforms}, then Proposition \ref{mainlemma} 
yields $\dim\mathcal{N}(\beta)\geq 2n-2p-4>0$. But this is a contradiction
since from \eqref{condition} we have that $\mathcal{N}(\beta)=0$.\qed  

\begin{proposition}\label{diag}\po Let the bilinear form 
$\beta\colon V^{2n}\times V^{2n}\to W^{p,p}$, $s\leq n$, be flat.
Assume that the vector subspace $\Sal(\beta)$ is nondegenerate and 
that \eqref{product} holds. For $p\geq 4$ assume further that there is 
no non-trivial $J$-invariant vector subspace $V_1\subset V^{2n}$ 
such that the subspace $\Sal(\beta|_{V_1\times V_1})$ is degenerate 
and $\dim\Sal(\beta|_{V_1\times V_1})\leq\dim V_1-2$.
Then $s=n$ and there is an orthogonal basis $\{X_i,JX_i\}_{1\leq i\leq n}$ 
of $V^{2n}$ such that:
\begin{itemize}
\item[(i)] $\beta(Y_i,Y_j)=0\;\mbox{if}\; i\neq j
\;\mbox{and}\;Y_k\in\spa\{X_k,JX_k\}$\;\mbox{for}\;k=i,j. 
\item[(ii)] The vectors 
$\{\beta(X_j,X_j),\beta(X_j,JX_j)\}_{1\leq j\leq n}$ 
form an orthonormal basis of $\Sal(\beta)$.
\end{itemize}
\end{proposition}

\proof It follows from  Proposition $15$ in \cite{CD}.
\vspace{2ex}\qed

\noindent{\em Proof of Theorem \ref{three}:} 
Theorem \ref{one} gives that $p=n-2$. In an open neighborhood 
$U$ of $x_0$ in $M^{2n}$ there is a  complex vector subbundle 
$\bar{V}\subset TM$ such that $\bar{V}(x_0)=V^{2m}$ and 
$K_M(S,JS)>0$ for any $0\neq S\in\bar{V}$. At any point of $U$ 
the vector subspace $\Sal(\beta)$ is nondegenerate. In fact, 
if otherwise then by Proposition \ref{umbilical} there 
is a point $y\in U$ and a complex vector subspace 
$P^{2\ell}\subset T_yM$ with $\ell\geq 2$ such that the 
sectional curvature satisfies $K_M(S,JS)\leq 0$ for any 
$0\neq S\in P^{2\ell}$, in contradiction with our assumption.

By Proposition \ref{diag}, there is at any $y\in U$ an
orthogonal basis $\{X_j,JX_j\}_{1\leq j\leq n}$ of $T_yM$ 
such that both parts hold. By part $(ii)$ the vectors
$(\xi_j,0)=\beta(X_j,X_j)\in N_FM(y)$, $1\leq j\leq n$, 
are orthonormal. Then the argument used for the proof 
of Lemma $18$ in \cite{CD} gives that $F|_U$ has flat normal 
bundle, that $\rank A_{\xi_j}$=2 for $1\leq j\leq n$ and 
that the normal vector fields $\xi_1,\ldots,\xi_n$ are 
smooth on connected components of an open dense subset of $U$.
Moreover, we obtain from the Codazzi equation and the use 
of the de Rham theorem that $M^{2n}$ is locally a Riemannian 
product of surfaces $M_1^2\times\cdots\times M_n^2$.

Having that the codimension is $n=p+2$ and that $\a^F(Y_i,Y_j)=0$
if $Y_i\in (E_i)$ and $Y_j\in (E_j)$, $i\neq j$, then by 
Theorem $8.7$  in \cite{DT} there are isometric immersions 
$g_1\colon M_1^2\to\L^3$ and
$g_j\colon M_j^2\to\R^3$, $2\leq j\leq n$, such that
$$
F(x_1,\ldots,x_n)=(g_1(x_1),g_2(x_2),\ldots,g_n(x_n)).
$$
Since $F(M)\subset\mathbb{V}^{3n-1}\subset\L^{3n}$ then
$\<F,F\>=0$. Hence $\<g_j{}_*X_j,g_j\>=\<g_*X_j,g_j\>=0$ 
and thus $\|g_j\|=r_j$ with $-r_1^2+\sum_{j=2}^nr_j^2=0$.
This gives that $F(U)\subset
\Hy^2_{c_1}\times\Sp^2_{c_2}\times\cdots\times\Sp^n_{c_n}$
where $1/c_i=r_i^2$ and, by continuity, this also holds 
for $F(M)$. \vspace{2ex}\qed

\subsection*{Acknowledgments}

Luis J. Al\'{\i}as and Marcos Dajczer are  partially supported 
by the grant PID2021-124157NB-I00 funded by 
MCIN/AEI/10.13039/501100011033/ `ERDF A way of making Europe',
Spain, and are also supported by Comunidad Aut\'{o}noma de la Regi\'{o}n
de Murcia, Spain, within the framework of the Regional Programme
in Promotion of the Scientific and Technical Research (Action Plan 2022),
by Fundaci\'{o}n S\'{e}neca, Regional Agency of Science and Technology,
REF, 21899/PI/22.

\subsection*{Data availability} 
Data sharing is not applicable to this article as no datasets were generated or analyzed during the current study.

\subsection*{Conflict of interest}
No potential conflict of interest was reported by the authors.

\noindent  Luis J. Al\'{i}as\\
Departamento de Matem\'{a}ticas\\
Universidad de Murcia,\\
E-30100 Espinardo, Murcia -- Spain\\
e-mail: ljalias@um.es
\bigskip

\noindent Sergio J. Chion Aguirre\\
Pontificia Universidad Cat\'{o}lica del Per\'{u},\\
CENTRUM Cat\'{o}lica Graduate Business School,\\
Lima -- Per\'{u},\\
e-mail: sjchiona@pucp.edu.pe
\bigskip

\noindent Marcos Dajczer\\
Departamento de Matem\'{a}ticas\\
Universidad de Murcia,\\
E-30100 Espinardo, Murcia -- Spain\\
e-mail: marcos@impa.br
\end{document}